\documentclass[12pt]{amsart}
\usepackage{graphicx, amssymb, amsthm, amsfonts, latexsym, epsfig,
amsmath, amscd, amsbsy}
\usepackage[usenames]{color}

\theoremstyle{definition}
\newtheorem{theorem}{Theorem}[section]
\newtheorem{lem}[theorem]{Lemma}

\newtheorem{pro}[theorem]{Proposition}
\newtheorem{df}[theorem]{Definition}
\newtheorem{rem}[theorem]{Remark}

\numberwithin{equation}{section}

\newenvironment{proofof}[1]{\medskip\noindent{\emph{Proof of #1.}}}{
\hfill\qed\\ }

\def\R{{\mathbb R}}
\def\N{{\mathbb N}}
\def\Z{{\mathbb Z}}

\def\C0{{\mathfrak C}}
\def\Int{\mathop\mathrm{Int}}
\def\Cl{\mathop\mathrm{Cl}}

\def\ILq{\underleftarrow\lim([u,c_1],q_a)}
\def\ILs{\underleftarrow\lim([0,1],T_{s})}
\def\ILcores{\underleftarrow\lim([c_2, c_1],T_{s})}

\def\htop{h_{top}}

\def\chain{{\mathcal C}}

\def\ie{{\em i.e.,}\ }

\def\eps{\varepsilon}
\def\diam{\mbox{diam}}

\def\mesh{\mbox{mesh}}

\def\orb{\mbox{orb}}

\def\lap{\mathfrak l}

\begin{document}

\title[Entropy of Unimodal Inverse Limits]{Entropy of homeomorphisms on
unimodal inverse limit spaces}
\author{H.~Bruin and S.~\v{S}timac}
\thanks{HB was supported by EPSRC grant EP/F037112/1.
S\v{S} was supported in part by the MZOS Grant 037-0372791-2802 of the
Republic of Croatia.}
%\date{\today}

\subjclass[2000]{54H20, 37B45, 37E05}
\keywords{entropy, inverse limit space, tent map, logistic map}

\begin{abstract}
We prove that every self-homeomorphism $h : K_s \to K_s$ on the inverse
limit space $K_s$ of the tent map $T_s$
with slope $s \in (\sqrt 2, 2]$ has topological entropy $\htop(h) = |R|
\log s$, where $R \in \Z$ is such
that $h$ and $\sigma^R$ are isotopic.
Conclusions on the possible values of the entropy of homeomorphisms of
the inverse limit space of a (renormalizable) quadratic map are drawn as
well.
\end{abstract}

\maketitle

\baselineskip=18pt

\section{Introduction}\label{sec:intro}

Since Williams' work on expanding attractors \cite{W}, inverse limit
spaces play a key role in constructing
and analyzing examples in dynamical systems. For a continuous map $f : X
\to X$ of a metric space $X$, its
shift homeomorphism $\sigma_f : X_f \to X_f$, where $X_f := {\displaystyle
\lim_{\longleftarrow}(X, f)}$,
is a homeomorphism of the inverse limit space $X_f$, and it is the
dynamically minimal extension of $f$ to
a homeomorphism. Although the space $X_f$ is defined as a subspace of
$X^{\N}$, Whitney's Embedding Theorem
allows one to embed the inverse limit in a manifold $M$ (of dimension $3$
if $X \subset \R$), and
$\sigma_f : X_f \to X_f$ can be extended to a homeomorphism of $M$ for
which $X_f$ is a global attractor.

Barge and Martin \cite{BM} described a construction to embed the inverse limit of any interval endomorphism
as a global attractor of a planar homeomorphism and their construction readily generalizes to graphs other
than the interval. Bruin \cite{B2} showed that for unimodal maps, this planar homeomorphism can be taken
Lipschitz continuous. Very recently, Boyland, de Carvalho and Hall have developed a $C^0$ parametrized version
\cite{BCH} of the Barge-Martin construction. In the same paper, they apply their results to the one of the
most studied family of unimodal maps, the tent family $T_s : [0, 1] \to [0, 1]$ with slope $\pm s$,
$s \in [1, 2]$, defined as $T_s(x) = \min\{sx, s(1-x)\}$. For this family it has been proven that for
$1 \le s < t  \le 2$, the inverse limits $K_s := \underleftarrow\lim([0, 1],T_s)$ and $K_t$ are not
homeomorphic, see \cite{Kail2, Stim, RS, BBS}. By Barge and Martin \cite{BM}, any map of the interval can be
extended to the disk, and their construction gives a family of extensions of any unimodal family. Thus
Boyland, de Carvalho, and Hall provide a continuously varying family $\phi_s : {\mathbb D} \to {\mathbb D}$
of homeomorphisms of the disk ${\mathbb D}$, with $\phi_s$ having an embedded copy of $K_s$ as a global
attracting set; the dynamics of $\phi_s$ restricted to this attractor is conjugate to the shift homeomorphism $\sigma := \sigma_{T_s}:K_s \to K_s$.
Moreover, by \cite{BCH}, these attractors vary continuously in the Hausdorff topology, although no two are
homeomorphic.

In this paper, we add to the description of the group of homeomorphisms $h : K_s \to K_s$ for arbitrary $s$.
In \cite{BJKK, BKRS} and \cite{BS} it is proven that for $s \in (\sqrt 2, 2]$ and a
homeomorphism $h : K_s \to K_s$, there exists $R \in \Z$ such that $h$ is isotopic to $\sigma^R$.
Here we study the topological entropy $\htop(h)$ and show:

\begin{theorem}\label{thm:entropy}
If $s \in (\sqrt2, 2]$, $h : K_s \to K_s$ is a homeomorphism and $R \in \Z$ is such that $h$ is isotopic to
$\sigma^R$, then the topological entropy $\htop(h) = |R| \log s$.
\end{theorem}

Isotopy does in general not preserve topological entropy, as is easy to see in dimension $\ge 2$. The
one-dimensional situation of Theorem~\ref{thm:entropy} is more rigid, and one is inclined to compare it with
skew-products where all fiber maps are monotone interval maps.
Indeed \cite{BS}, the rigidity of isotopy implies that
$h$ can differ from $\sigma^R$ only at {\em non-folding points,}
(see Definition~\ref{def:folding_point} below).
Using Bowen's entropy formula and the notion of
sequence entropy, Kolyada, Misiurewicz and Snoha \cite{KMS,KS} proved that for such skew-products the entropy
is the same as the entropy of the base map. In our case, we cannot easily identify the base map, and at best
it would be defined on a non-compact union of Cantor sets. Therefore it is unclear if results from the
skew-product setting can be applied and a proof must be given. In this paper we give our own proof.

We call the fact that the homeomorphisms on a space $X$ can only taking
specific values {\em entropy rigidity}. Clearly the circle has entropy rigidity:
zero is the only possible entropy.
A related example is the pseudo-arc which can be written as inverse limit
space of an interval map in various ways. Mouron \cite{M}
proved that zero and infinity are the only possible values of entropy
for the corresponding shift homeomorphism.

The entropy rigidity result of Theorem~\ref{thm:entropy} can be extended
to the inverse limit spaces of
quadratic maps $q_a(x) = 1-ax^2$. Each $q_a$, $a \in [0, 2]$, with
positive topological entropy, is
semi-conjugate to a tent map $T_s$ with $\log s = \htop(q_a)$ \cite{MT}
and  the semi-conjugacy collapses (pre)periodic intervals to points.
If $q_a$ has no non-trivial periodic intervals (\ie $q_a$ is not
renormalizable), then Theorem~\ref{thm:entropy} applies to its inverse
limit space.
Otherwise, the first return map to such a periodic interval is a new
unimodal map, called the {\em renormalization} of the previous; thus we
get a (possibly infinite)
sequence of nested cycles of periodic intervals, with periods $(p_i)_{i
\ge 0}$, where $p_i$ divides
$p_{i+1}$ and $p_0 = 1$. The effect of renormalization on the structure of
the inverse limit space is
well-understood \cite{BaDi3}: it produces proper subcontinua that are
$p_i$-periodic under the shift
homeomorphism and homeomorphic with the inverse limit space of the
renormalized map. Here we show that even in this case an analogous of 
Theorem~\ref{thm:entropy} holds:

\begin{theorem}\label{thm:entropy_quadratic}
Assume that $q_a$ is a quadratic map with positive topological entropy and $\log s = \htop(q_a)$.
If $h$ is a homeomorphism on the inverse limit space of $q_a$, then the topological entropy 
$\htop(h) = |R|\log s$ where $R \in \Z$ is such that $h$ is isotopic to $\sigma^R$.
\end{theorem}

As unimodal inverse limit spaces mimic to some extent the structure of H\'enon attractors,
Theorems~\ref{thm:entropy} and \ref{thm:entropy_quadratic} can be seen as a step towards an entropy rigidity
result for H\'enon attractors. Recall that the H\'enon map is $H_{a,b}(x, y) = (1 + y - ax^2, bx)$, and for
parameters $a \in (0,2)$ and small $b$, there is a bounded open forward invariant set $U$ such that
$\cap_{n \ge 0} H^n(U)$ is a continuum, called the {\em (global) attractor}. Barge and Holte \cite{BaH} proved
that for parameter values $a$ where $q_a$ has an attracting periodic orbit, and small values of $b$, $H_{a,b}$
restricted to the global attractor is topologically conjugate to shift homeomorphisms on inverse limit of an
interval with bonding map $q_a$. Hence Theorem~\ref{thm:entropy_quadratic} applies in this situation,
but due to a result by Barge \cite{barge}, this may be the only case where H\'enon attractors
are homeomorphic to unimodal inverse limit spaces with a single bonding map.

Despite this, it would be interesting to see if general H\'enon (or Lozi) attractors exhibit rigidity of
entropy. Specifically, what is the class of H\'enon maps $H_{a,b}$ for which entropy of every
self-homeomorphism of the corresponding attractor is an integer multiple of $h_{top}(H_{a,b})$?

Let us mention a final example of entropy rigidity pertaining to certain skew-products of tent maps, see
\cite{HRS}. In this case, the attractor is a projection of the inverse limit space of a tent map, and the
entropy rigidity is as described in Theorem~\ref{thm:entropy}.
\\[4mm]
{\em Acknowledgement:} We would like to thank the referees for their very useful
comments on this paper.

\section{Preliminaries}\label{sec:prelim}

Let $\N = \{ 1,2,3,\dots\}$ be the set of natural numbers and $\N_0 = \N
\cup \{ 0 \}$. We consider the
family of tent maps $T_s:[0,1] \to [0,1]$ with slope $\pm s$, $s \in [1,
2]$, defined as
$T_s(x) = \min\{sx, s(1-x)\}$. The critical or turning point is $c :=
1/2$. We write $c_k := T_s^k(c)$, so
in particular $c_1 = s/2$ and $c_2 = s(1-s/2)$. The closed $T_s$-invariant
interval $[c_2, c_1] = [s-s^2/2, s/2]$ is called
the {\em core}.

The inverse limit space $K_s := \ILs$ is the collection of all backward
orbits
\[
\{ x = (\dots, x_{-2}, x_{-1}, x_0) :  T_s(x_{-i-1}) = x_{-i} \in [0,s/2]
\textrm{ for all } i \in \N_0 \},
\]
equipped with metric $d(x,y) = \sum_{n \le 0} 2^n |x_n - y_n|$ and {\em
induced} $($or {\em shift$)$
homeo\-morphism}
\[
\sigma_{T_s}(x) :=  \sigma(\dots, x_{-2}, x_{-1}, x_0) = (\dots, x_{-2}, x_{-1}, x_0, T_s(x_0)).
\]
Let $\pi_k: \ILs \to [0, c_1]$, $\pi_k(x) = x_{-k}$ be the $k$-th
projection map.
We fix
$s \in (\sqrt{2}, 2]$; for these parameters $T_s$ is not renormalizable and
$\underleftarrow\lim([c_2, c_1],T_s)$ is indecomposable. For any point $x \in K_s$,
the composant of $x$ in $K_s$ is the union of all proper subcontinua of $K_s$
containing $x$, and the arc-component of $x$ in $K_s$ is the union of all arcs in $K_s$
containing $x$.

We review some of the main tools introduced in \cite{BBS} and which are necessary here as well.
Since $0 \in [0, c_1]$ is
a fixed point of $T_s$, the endpoint $\bar 0 := (\dots, 0,0,0)$ is
contained in  $\ILs$.
The arc-component of $K_s$ which contains $\bar 0$ will be denoted as $C$;
it is a ray converging to, but
disjoint from the inverse limit of the core $\underleftarrow\lim([c_2,
c_1],T_s)$.

We define {\em $p$-points} as those points
$x = (\dots, x_{-2}, x_{-1}, x_0) \in K_s$  such that $x_{-p-k} = c$ for some
$k \in \N_0$. The number $L_p(x) := k$ is called the {\em $p$-level} of
$x$.  By convention, the endpoint $\bar 0$ of $C$ is also a $p$-point and
$L_p(\bar 0) := \infty$, for every $p$.

The ordered set $\{ z^1, z^2, \dots , z^n, \dots \}$ of all $p$-points of the arc-component $C$ is denoted
by $E_p$, and the ordered set of all $p$-points of the
arc-component $C$ of $p$-level $k$ by $E_{p,k}$.
The sequence of $p$-levels $( L_p(x) )_{x \in E_p}$ is called the
{\em folding pattern} of $C$, because it indicates the way how $C$ ``folds''
back and forth when it compactifies on the core $\ILcores$.
We denote the folding pattern as
$$
FP_p(C) = L_p(z^1), L_p(z^2), \dots , L_p(z^n), \dots
$$
Let $q \in \N$, $q > p$, and
$E_q = \{ y^1, y^2, \dots , y^n, \dots \}$. Since $\sigma^{q-p}$ is an
order-preserving homeomorphism of $C$,
it is easy to see that $\sigma^{q-p}(z^i) = y^i$ for every $i \in \N$, and
$L_p(z^i) = L_q(y^i)$. Therefore,
the folding pattern of $C$ does not depend on $p$.
In fact, the graph of $T^p_s|_{[0,c_1]}$ oscillates in the same way
as the initial arc of $C$, and therefore $FP_p(C)$
starts as $\infty \ 0 \ 1 \ 0 \ 2 \ 0 \ 1 \ \dots$ for every
slope $s > 1$ and $p \in \N_0$.

More generally, given an arc $A \subset K_s$ with successive $p$-points
$x^0, \dots , x^n$, the sequence of their $p$-levels is denoted as
$$
FP_p(A) := L_p(x^0), \dots , L_p(x^n).
$$
Note that every arc of $C$ has only finitely many $p$-points, but an
arc $A$ of the core of $K_s$ can have infinitely many $p$-points. In this
case, if $(u^i)_{i \in {\mathcal I}}$
is the set of $p$-points of $A$, then $FP_p(A) = (L_p(u^i))_{i \in
{\mathcal I}}$, for some countable index
set $\mathcal I$ (not necessarily of the same ordinal type as $\N$ or
$\Z$).

In \cite{BBS} it is shown that the asymptotic structure of the
folding pattern of $C$ is a topological invariant, but it takes a long argument
to reach to that conclusion.
For instance, $p$-points have no topological characterization and a homeomorphism will in general not map $p$-points to $q$-points for
any integers $p, q$.
However, among the $p$-points there are special ones, which we call
{\em salient}, which are center points of symmetries in $C$.
Homeomorphisms preserve these symmetries to such an extent that it is possible
to prove that salient points map close to salient points.

\begin{df}\label{df:salient}
We call a $p$-point $y \in C$ \emph{salient} if $0 \le L_p(x) < L_p(y)$ for every $p$-point
$x \in (\bar 0 , y)$. Let $(s_p^i)_{i \in \N}$ be the sequence of all salient $p$-points of $C$,
ordered such that $s_p^i \in (\bar 0 , s_p^{i+1})$ for all $i \ge 1$.
\end{df}

Since by definition $L_p(s_p^i) > 0$, for all $i \ge 1$, we have $L_p(s_p^1) = 1$. Also, since
$s_p^i = \sigma^{i-1}(s_p^1)$, we have $L_p(s_p^i) = i$, for every $i \in \N$. Therefore, for every
$p$-point $x$ of $K_s$ with $L_p(x) \ne 0$, there exists a unique salient $p$-point $s_p^k$ such that
$L_p(x) = L_p(s_p^k) = k$. Also, for every $k \in \N$, among all $p$-points $E_{p, k}$ of $C$ with
$p$-level $k$ there exists precisely one $p$-point $s_p^k$ which is salient and has $p$-level $k$. Note
that the salient $p$-points depend on $p$: if $p \ge q$, then the salient $p$-point $s_p^i$ equals the salient $q$-point $s_q^{i+p-q}$.

\begin{df} \label{def:folding_point}
A {\em folding point} is any point $x$ in the core of $K_s$ such that no neighborhood of $x$ in the core
of $K_s$ is homeomorphic to the product of a Cantor set and an arc.
\end{df}

In \cite{BS} it was shown that
$x \in K_s$ is a folding point if and only if for some $p \in \N$ there is a sequence of $p$-points
$(x^k)_{k \in \N}$ such that $x^k \to x$ and $L_p(x^k) \to \infty$.

\section{Construction of chains $\chain_p$ and $\widetilde
\chain_{p+M}$}\label{sec:construction}

A continuum is {\em chainable} if for every $\eps > 0$, there is a cover $\{ \ell^1, \dots , \ell^n\}$ of
open sets (called {\em links}) of diameter $< \eps$ such that $\ell^i \cap \ell^j \neq \emptyset$ if and only
if $|i-j| \le 1$. Such a cover is called a {\em chain}. Clearly the interval $[0,s/2]$ is chainable.
Throughout, we will use sequence of chains $\chain_p$ of $\ILs$ satisfying the following properties:
\begin{enumerate}
\item there is a chain $\{ I^1_p, I^2_p, \dots , I^n_p \}$ of $[0,s/2]$ such that
$\ell^j_p := \pi_p^{-1}(I^j_p)$ are the links of $\chain_p$;
\item each point $x \in \cup_{i=0}^p T_s^{-i}(c)$ is a boundary point of some link $I^j_p$;
\item for each $i$ there is $j$ such that $T_s(I^i_{p+1}) \subset I^j_p$.
\end{enumerate}
If $\max_j |I^j_p| < \eps s^{-p}/2$ then
$\mesh(\chain_p) := \max\{ \diam(\ell_p) :\ell_p \in \chain_p\} < \eps$,
which shows that $\ILs$ is indeed chainable. Condition (3) ensures that $\chain_{p+1}$ \emph{refines}
$\chain_p$ (written $\chain_{p+1} \preceq \chain_p$).

Note that all $p$-point $E_{p, k}$ of $p$-level $k$ belong to the same link of $\chain_p$. (This follows
by property (1) of $\chain_p$, because $L_p(x) = L_p(y)$ implies $\pi_p(x) = \pi_p(y)$.) Therefore, every
link of $\chain_p$ which contains a $p$-point of $p$-level $k$, contains also the salient $p$-point $s_p^k$.

Let $h:\ILs \to \ILs$ be a homeomorphism. Take $q, p \in \N_0$ such that
\[
h(\chain_{q}) \preceq \chain_{p}.
\]
This means that it becomes natural to measure the $p$-level of $p$-points
close to the images $h(s_q^i)$ of salient $q$-points.
Let us denote by $\ell^x_p$ a link of $\chain_p$ which contains the point $x$. From \cite{Stim} and \cite{BBS}
(for the finite and infinite critical orbit case, respectively) we have the following proposition:

\begin{pro}\label{pro:levels}
There exists $M \in \Z$ such that the following holds:
\begin{itemize}
\item[(i)] Let $k \in \N$ and let $s_q^k$ be a salient $q$-point.
Then $u := h(s_q^k) \in \ell_p^{s_p^{k+M}}$ and the arc component $A_u \subset \ell_p^{s_p^{k+M}}$
containing $u$, also contains the salient $p$-point $s_p^{k+M}$.
\item[(ii)] Let $k \in \N$ and let $x'$ be a $q$-point with $L_q(x') = k$.
Then $u := h(x') \in \ell_p^{s_p^{k+M}}$ and the arc component $A_u \subset \ell_p^{s_p^{k+M}}$ containing $u$,
also contains a $p$-point $x$ such that $L_p(x) = k + M$.
Moreover, the number of $q$-points in $[s_q^i, s_q^{i+1}]$ with $q$-level $k$, $k \le i$, is the same as
the number of $p$-points in $[s_p^{M+i}, s_p^{M+i+1}]$ with $p$-level $M+k$, for every $i \in \N$.
\end{itemize}
\end{pro}

{}From \cite{BJKK, BKRS} and \cite{BS} we can derive

\begin{pro}\label{prop:independent} The integer $R = M + p - q$ does not depend on $M, p, q$, and is such
that $h$ and $\sigma^R$ are isotopic.
\end{pro}

Let us write $x \approx_p y$ if $x$ and $y$ belong to the same arc component of the same link of the chain
$\chain_p$. Using this notation we can write Proposition~\ref{pro:levels} in the following way:
$h(x) \approx_p \sigma^R(x)$ for every $x \in E_q$. A fortiori, it was shown in \cite{BS} that $\sigma^{-R} \circ h$
is the identity on the set of folding points of $\ILcores$. At non-folding points, there is flexibility that
a priori might allow a difference in entropies of $h$ and $\sigma^R$. For an  arbitrary $m \in \N$, the
$h^m$-images of links of $\chain_q$
are hard to control;
for a link $\ell \in \chain_q$, images $h^m(\ell)$ and $\sigma^{Rm}(\ell)$ need not even be
close to ``parallel'' in $\chain_p$, where ``parallel'' in $\chain_p$ means
$\pi_g(h^m(\ell)) = \pi_g(\sigma^{Rm}(\ell))$ for $g < p$.
To overcome this problem, we need to introduce new chains $\widetilde \chain_{p+M}$ with very elongated links
(being the concatenation of links of $\chain_{p+M}$)
such that $h^m$-image of $\chain_q$
can be shown to be ``parallel'' to $\widetilde \chain_{q+Rm} = \widetilde \chain_{p+M}$ (we define the notion
`basically collinear' in Remark~\ref{rem:collinear2} to make this precise).
The rest of this section is devoted to the properties of these additional chains.

Take $p < q \in \N_0$ such that $h^m(\chain_q) \prec \chain_p$ and $\sigma^{Rm}(\chain_q) \prec \chain_p$
and let $M = Rm + q - p$. For each $j \ge 1$, $E_{q,j}$ is contained in a single link $\ell_q \in \chain_q$
and by Proposition~\ref{pro:levels}, for $\ell_p \supseteq h^m(\ell_q)$, every point of $h^m(E_{q,j})$ is
contained in an arc component of $\ell_p$ which contains a $p$-point of $E_{p,M+j} = E_{p+M,j}$.

By using $h^{-1}$ instead of $h$ if necessary, we can assume that $R \ge 0$. Then $M > 0$ and
$\chain_{p+M} \prec \chain_p$. Let $\widetilde \chain_{p+M}$ be a chain with
$\chain_{p+M} \preceq \widetilde \chain_{p+M} \preceq \chain_p$ satisfying the following property: Each link
$\widetilde \ell_{p+M} \in \widetilde \chain_{p+M}$ is much wider in the `tangential' direction than in the
`transversal' direction; in fact, the `tangential' length of $\widetilde \ell_{p+M}$ is the same as the
`tangential' length of $\ell_p \in \chain_p$, and the `transversal' length of
$\widetilde \ell_{p+M} = \widetilde \ell_{q+Rm}$ (recall that $p + M = q + Rm$) is the same as the
`transversal' length of $\ell_{q+Rm} \in \chain_{q+Rm}$, see Figure~\ref{fig:tildeCpM}.

\begin{figure}[ht]
\unitlength=8mm
\begin{picture}(14,6.5)(2,0)
{\color{blue}\put(7, 3){\oval(12,6)}}\put(0.4,5.5){$\ell_p$}
\put(1,2){\line(1,0){12}} \put(1,2.1){\line(1,0){12}}
\put(1,2.3){\line(1,0){12}} \put(1,2.4){\line(1,0){12}}
\put(1,3){\line(1,0){12}} \put(1,3.1){\line(1,0){12}}
\put(1,3.3){\line(1,0){12}} \put(1,3.4){\line(1,0){12}}
\put(10, 5){\oval(0.5,0.5)[l]} \put(10,4.75){\line(1,0){3}}
\put(10,5.25){\line(1,0){3}} \put(10, 5){\oval(0.3,0.3)[l]}
\put(10,4.85){\line(1,0){3}} \put(10,5.15){\line(1,0){3}}
\put(11, 1){\oval(0.5,0.5)[l]} \put(11,0.75){\line(1,0){2}}
\put(11,1.25){\line(1,0){2}} \put(11, 1){\oval(0.3,0.3)[l]}
\put(11,0.85){\line(1,0){2}} \put(11,1.15){\line(1,0){2}}
\put(5.5,5){$\ell^j_{p+M} \in \chain_{p+M}$}
\put(5.5,4.7){\vector(-2,-3){0.8}}\put(5.7,4.7){\vector(0,-1){1.2}}
\put(6,4.7){\vector(2,-1){2.4}}
\put(4.2,4.1){$\dots$}\put(6,4.1){$\dots$}\put(7.6,4.1){$\dots$}
\thicklines
{\color{red}
\put(0.95,1.8){\dashbox{0.2}(12.1,0.8)}
\put(0.95,2.8){\dashbox{0.2}(12.1,0.8)}
\put(9.45,4.5){\dashbox{0.2}(3.6,1)}
\put(10.45,0.5){\dashbox{0.2}(2.6,1)}}
\thinlines
\put(1.2, 2.2){\oval(1,0.6)} \put(2.1, 2.2){\oval(1,0.6)} \put(3.0,
2.2){\oval(1,0.6)} \put(3.9,2.2){\oval(1,0.6)} \put(4.8,
2.2){\oval(1,0.6)} \put(5.7,2.2){\oval(1,0.6)} \put(6.6,
2.2){\oval(1,0.6)} \put(7.5,2.2){\oval(1,0.6)} \put(8.4,
2.2){\oval(1,0.6)} \put(9.3,2.2){\oval(1,0.6)} \put(10.2,
2.2){\oval(1,0.6)} \put(11.1,2.2){\oval(1,0.6)} \put(12.0,
2.2){\oval(1,0.6)} \put(12.9, 2.2){\oval(1,0.6)}
\put(1.2, 3.2){\oval(1,0.6)} \put(2.1, 3.2){\oval(1,0.6)} \put(3.0,
3.2){\oval(1,0.6)} \put(3.9,3.2){\oval(1,0.6)} \put(4.8,
3.2){\oval(1,0.6)} \put(5.7,3.2){\oval(1,0.6)} \put(6.6,
3.2){\oval(1,0.6)} \put(7.5,3.2){\oval(1,0.6)} \put(8.4,
3.2){\oval(1,0.6)} \put(9.3,3.2){\oval(1,0.6)} \put(10.2,
3.2){\oval(1,0.6)} \put(11.1,3.2){\oval(1,0.6)} \put(12.0,
3.2){\oval(1,0.6)} \put(12.9, 3.2){\oval(1,0.6)}
\put(10.1, 5.0){\oval(1,0.7)} \put(11,5.0){\oval(1,0.7)} \put(11.9,
5.0){\oval(1,0.7)} \put(12.8, 5.0){\oval(1,0.7)}
\put(11.0,1.0){\oval(1,0.7)} \put(11.9, 1.0){\oval(1,0.7)}
\put(12.8, 1.0){\oval(1,0.7)}
\put(14.8,5){\vector(-1,0){1.7}}
\put(14.8,4.9){\vector(-1,-1){1.7}}
\put(15,4.8){$\tilde\ell^j_{p+M} \in \widetilde\chain_{p+M}$}
\put(14.8,4.6){\vector(-3,-4){1.7}} \put(14.8,4.5){\vector(-1,-2){1.7}}
\end{picture}
\caption{Impression of the links $\ell_p$ (blue),
$\tilde\ell^j_{p+M} \in \widetilde\chain_{p+M}$
(red dashed lines) and $\ell^j_{p+M} \in \chain_{p+M}$ (black).}
\label{fig:tildeCpM}
\end{figure}

To further understand the chains $\chain_q$, $\chain_{q+Rm}$, $\widetilde \chain_{q+Rm}$ and $\chain_p$,
note that for every link $\widetilde \ell_{q+Rm} \in \widetilde \chain_{q+Rm}$, there is a link
$\ell_p \in \chain_p$ such that every arc component $G \in \widetilde \chain_{q+Rm}$ is an arc component of
$\ell_p$ as well. In other words, if $G_a$ and $G_b$ are two arc components of $\ell_p$, then they can belong
to different links of $\widetilde\chain_{q+Rm}$, say $G_a \subset \widetilde\ell^a_{q+Rm}$ and
$G_b \subset \widetilde\ell^b_{q+Rm}$, but if arc components $G_a \subset \ell_p$ and
$\tilde G_a \subset \widetilde\ell^a_{q+Rm}$ are such that $G_a \cap \tilde G_a \ne \emptyset$, then
$G_a = \tilde G_a$. Therefore, all points of $h^m(E_{q,j})$ are contained in the link
$\widetilde\ell_{q+Rm} \in \widetilde \chain_{q+Rm}$ which contains $E_{p+M,j} = E_{q+Rm,j}$. Moreover, by
construction of $\widetilde \chain_{q+Rm}$, we have that
$h^m(\chain_q) \prec \widetilde \chain_{p+M} = \widetilde \chain_{q+Rm}$.

Let $\chain_q = (\ell_q^i)_{i=1}^{n_q}$ and
$\widetilde \chain_{q+Rm} = (\tilde \ell_{q+Rm}^i)_{i=1}^{n_{q+Rm}}$. Let $A_q = [\bar 0, s_q^1]$ and
$A_{q+mR} = [\bar 0, s_{q+mR}^1]$, where $s_q^1$ and $s_{q+mR}^1$ are the first salient $q$-point and
salient $(q+mR)$-point respectively. Then $\pi_q(A_q) = [0, c_1] = \pi_{q+mR}(A_{q+mR})$. Also
$\bar 0 \in \ell_q^1$, $\bar 0 \in \tilde \ell_{q+Rm}^1$, $s_q^1 \in \ell_q^{n_q}$ and
$s_{q+mR}^1 \in \tilde \ell_{q+Rm}^{n_{q+Rm}}$. Therefore, $\chain_{q+Rm}$ and $\widetilde \chain_{q+Rm}$
coil through $\chain_q$ in the same way with the same number of laps as $T^{Rm}_s$ maps $[0, c_1]$ to itself.
Also, $\sigma^{Rm}(\chain_q)$ is {\em collinear} with $\widetilde \chain_{q+Rm}$, \ie $\sigma^{Rm}(\chain_q)$
goes straight through $\widetilde \chain_{q+Rm}$. More precisely, for every $t \in [c_2, c_1]$ there exists a
unique $j$ such that either $\pi_{q+Rm}^{-1}(t) \subset \sigma^{Rm}(\ell_q^j)$, or
$\pi_{q+Rm}^{-1}(t) \subset \sigma^{Rm}(\ell_q^j) \cap \sigma^{Rm}(\ell_q^{j+1})$.

\begin{rem}\label{rem:collinear1}
Note that if $a, b \in \N$ are such that all of $h^m(\ell_q^i)$, $a \le i \le b$, are mapped into the same
link of $\widetilde \chain_{q+Rm}$, \ie if $\cup_{i=a}^bh^m(\ell_q^i) \subseteq \tilde\ell_{q+Rm}^j$ for
some $j$, it is possible that there exists $t \in [c_2, c_1]$ such that
$\pi_{q+Rm}^{-1}(t) \cap h^m(\ell_q^i) \ne \emptyset$ for all $a \le i \le b$ even if $b - a \ge 2$.
To avoid this, let us consider the chain $V = (v^j)_{j=1}^N$ such that
\begin{itemize}
\item $v^j = \cup_{k = a_j}^{a_{j+1}-1} \ell_q^k$ for an increasing sequence $(a_j)_{j=1}^{N+1}$;
\item $a_1 = 1$, $a_{N+1} = n_q + 1$ with $\ell_q^{a_{N+1}} = \emptyset$ by convention;
\item for every $j \in \{ 1, \dots , N \}$ we have $\cup_{k = a_j}^{a_{j+1}-1} h^m(\ell_q^k) \subseteq \tilde\ell_{q+Rm}^i$, for some $i \in \{1, \dots , n_{q+Rm} \}$;
\item  $h^m(\ell_q^{a_{j+1}}) \nsubseteq \tilde\ell_{q+Rm}^i$ except for $a_{N+1}$ (since $\emptyset = h^m(\ell_q^{a_{N+1}}) \subset \tilde\ell_{q+Rm}^i$).
\end{itemize}
Such a chain exists since $h^m(\chain_q) \prec \widetilde \chain_{q+Rm}$, and in fact, by construction we have
$h^m(V) \prec \widetilde \chain_{q+Rm}$.
\end{rem}

The following lemma is needed for the proof of the main theorem. In that proof we need the property that
$h^m(\chain_q)$ coils through $\chain_q$ in roughly the same way with the same number of laps as $T^{mR}_s$
maps $[0, c_1]$ into itself.

\begin{lem}\label{lem:collinear}
The image chain $h^m(\chain_q)$ is basically collinear with $\widetilde \chain_{q+Rm}$, \ie $h^m(\chain_q)$
goes straight through $\widetilde \chain_{q+Rm}$.
\end{lem}

\begin{rem}\label{rem:collinear2}
In the above lemma ``basically collinear'' means exactly that there exists a chain $V = (v^j)_{j=1}^N$ as in
Remark~\ref{rem:collinear1} for which every link $v^j$ is a union of consecutive links of $\chain_q$,
such that for every $t \in [c_2, c_1]$, there exists a unique $j$ with either
$\pi_{q+Rm}^{-1}(t) \subset h^m(v^j)$, or $\pi_{q+Rm}^{-1}(t) \subset h^m(v^j) \cap h^m(v^{j+1})$.
\end{rem}

\begin{proofof}{Lemma~\ref{lem:collinear}} Let $V$ be as in Remark~\ref{rem:collinear1}. By
Remarks~\ref{rem:collinear1} and \ref{rem:collinear2} we only need to prove that the number of links in $V$
is the same as the number of links in $\widetilde \chain_{q+Rm}$, \ie $N = n_{q+Rm}$. In this case, for every
link $v^j \in V$ there exist a unique link $\tilde \ell_{q+Rm}^j \in \widetilde \chain_{q+Rm}$ such that
$h^m(v^j) \subseteq \tilde \ell_{q+Rm}^j$, implying that $h^m(V)$ is collinear with $\widetilde \chain_{q+Rm}$.

Note firstly that by Theorem 1.1 of \cite{BS}, we have $h(x) = \sigma^R(x)$ for every folding point $x$ in
$K_s$. Therefore, if the critical point $c$ is dense in the core $[c_2, c_1]$, then $h \equiv \sigma^R$
on the inverse limit of the core $\underleftarrow\lim([c_2, c_1],T_s)$. Since $\sigma^{Rm}(\chain_q)$ is
collinear with $\widetilde \chain_{q+Rm}$, also $h^m(\chain_q)$ is collinear with $\widetilde \chain_{q+Rm}$.

(1) Let us assume that the critical point $c$ is not dense in the core $[c_2, c_1]$. The construction of
$\widetilde \chain_{p+M}$ and Proposition~\ref{pro:levels} imply that $h^m(x) \approx_{q+Rm} \sigma^{Rm}(x)$
for every $x \in E_q$, meaning that $h^m(\chain_q)|_C$ is basically collinear with
$\widetilde \chain_{q+Rm}|_C$, \ie that $h^m(\chain_q)|_C$ goes straight through
$\widetilde \chain_{q+Rm}|_C$.

(2) Let $v^j$ be any link of $V$, and let $A_j$ be an arc-component of $v^j$ which is contained in the core
of $K_s$, and does not contain any folding point. Since we want to prove that $N = n_{q+Rm}$, we need not
consider those arcs $B$ such that $h^m(B)$ is contained in only one link of $\widetilde \chain_{q+Rm}$.
Hence, by extending the arc $A_j$ to $A = [x,y] \supset A_j$ if necessary, we can assume without loss of
generality that $x$ and $y$ are $q$-points (and/or folding point), but $A$ contains no folding point in its
interior. Since $C$ is dense in the core, for every $\delta > 0$, there exists an arc $D_{\delta} \subset C$
such that $A$ and $D_{\delta}$ are $\delta$-close in the Hausdorff metric,
and that $FR_q(A) = FP_q(D_{\delta})$. (In the case that one of the boundary points is a folding point
$FR_q(\Int A) = FP_q(\Int D_{\delta})$.) By continuity of $h$, this implies that there exists an arc
$D \subset C$ such that $h^m(A)$ and $h^m(D)$ go through $\widetilde \chain_{q+Rm}$ in the same way.
Therefore, in this case it suffices to consider only arcs contained in $C$.

(3) Let us suppose that the arc-component $A_j$ of $v^j$ contains a folding point $z$ (if $A_j$ contains
more than one folding point, let $z$ and $z'$ be such that the arc $[z, z'] \subset A_j$ contains all folding
points of $A_j$). Let $A_{j-1} \subset v^{j-1}$ and $A_{j+1} \subset v^{j+1}$ be those arc-components for
which $A_j \cap A_{j-1} \ne \emptyset$ and $A_j \cap A_{j+1} \ne \emptyset$ (in the case that $A_j$
contains also a $q$-point, we have either $A_{j-1}, A_{j+1} \subset v^{j-1}$, or
$A_{j-1}, A_{j+1} \subset v^{j+1}$). Let us assume that $A_{j-1}$ and $A_{j+1}$ do not contain any folding
point. Since $h(u) = \sigma^R(u)$ for every folding point $u$ in $K_s$, it suffices to consider arcs
$A = [x, z] \supset A_{j-1}$ and $A' = [z', y] \supset A_{j+1}$ such that $A$ and $A'$ do not contain any
folding point in its interiors (if $A_j$ contains only one folding point, then $z = z'$). Now, as in (2), we
can find arcs $D, D' \subset C$ such that $h^m(A)$ and $h^m(D)$, as well as $h^m(A')$ and $h^m(D')$, go
through $\widetilde \chain_{q+Rm}$ in the same way.

(4) If $A_{j-1}$, or $A_{j+1}$, contains folding points, then either at least one of $A_{j-1} \setminus A_j$,
$A_{j+1} \setminus A_j$ is contained in the single link of $V$, or we can find $j_1 \le j-1$ (and/or
$j_2 \ge j+1$) such that $A_{j_1}$ (and/or $A_{j_2}$) does not contain any folding point. In the first
case, let us suppose that $A_{j-1} \setminus A_j$ is contained in a single link of $V$. Then
$h^m(A_{j-1} \setminus A_j)$ is contained in the single link of $\widetilde \chain_{q+Rm}$, implying that we
need not consider this arc (to prove that $N = n_{q+Rm}$, it suffices to consider those arcs $B$ such that
$h^m(B)$ goes through several link of $\widetilde \chain_{q+Rm}$). The same conclusion follows if
$A_{j+1} \setminus A_j$ is contained in a single link of $V$. In the second case we can proceed as in (3).

From the above steps we conclude that it is sufficient to consider only arcs contained in $C$. Since by (1),
$h^m(\chain_q)|_C$ is basically collinear with $\widetilde \chain_{q+Rm}|_C$, this implies that
$N = n_{q+Rm}$. This completes the proof.
\end{proofof}

\section{Proofs of Theorems~\ref{thm:entropy} and
\ref{thm:entropy_quadratic}} \label{sec:proof}

We will use Bowen's definition of $(n,\eps)$-separated sets to compute topological entropy of $h$.
For interval maps $T_s$, the  entropy is also the exponential growth-rate of the lap-number of $T^n$,
\ie $\htop(T_s) = \lim_{n \to \infty} \frac1n \log \lap(T_s^n)$, where
$$
\lap(T_s^n) := \#\{\text{maximal intervals on which $T^n_s$ is monotone}\}.
$$

\begin{proofof}{Theorem~\ref{thm:entropy}}
Fix homeomorphism $h$ and let $R$ be such that $h$ and $\sigma^R$ are isotopic. By using $h^{-1}$ instead of
$h$ if necessary, and noting that $\htop(h) = \htop(h^{-1})$, we can assume that $R \ge 0$.

We start with the upper bound. Take $\eta > 0$ arbitrary, and find $m \in \N$ such that the lap-number
$\lap(T_s^{Rm}) \le s^{m(R + \eta)}$ and also $(\log 2)/m < \eta$.

Fix $\eps > 0$ and let $p \in \N$ be such that $\mesh(\chain_p) < \eps/2$. Find $\eps_0, \eps_1 > 0$ such
that $\eps_1 < \eps_0/200 < \eps_0 < \mesh(\chain_p)$. Take $q \in \N$ such that $\mesh(\chain_q) < \eps_1$,
$h^m(\chain_q) \prec \chain_p$ and $\sigma^{Rm}(\chain_q) \prec \chain_p$.

Let $G = (g^j)_{j=1}^{N'}$ be a `chain' composed of `half-open links'
$g^j = (\cup_{k = a_j}^{a_{j+1}} \ell_q^k) \setminus \ell_q^{a_{j+1}}$ for an increasing sequence
$(a_j)_{j=1}^{N'}$ such that the following two conditions hold:
\begin{itemize}
\item[(a)] $a_1 = 1$, $a_{N'+1} = n_q+1$ and $\ell_q^{a_{N'+1}} =
\emptyset$ by convention,
\item[(b)] $\cup_{k = a_j}^{a_{j+1}-1} \ell_q^k \subseteq \ell_p^j$, for
every $j \in \{ 1, \dots  N' \}$.
\end{itemize}
Note that the sequence $(a_j)_{j=1}^{N'}$ is chosen so that
$$
\eps_0 < \diam(g^j) < \eps/2.
$$
Since $\diam(\ell^k_q) < \eps_1$ there must be at least $\eps_0/\eps_1 \ge 200$ links $\ell^k_q$ inside each
$g^j$, and $g^j$ is much wider in the `tangential' direction than in the `transversal' direction. Although
$G$ is not a chain since $g^i \cap g^j = \emptyset$ for every $i \ne j$, we have the intersection of closures
$\Cl(g^i) \cap \Cl(g^j) \ne \emptyset$ if and only if $|i-j| \le 1$. Clearly, $\bar 0 \in g^1$ and
$s_q^1 \in g^{N'}$, where $s_q^1$ is the first salient $q$-point. Note that $G$ is coarser than $V$ from
Remark~\ref{rem:collinear1}, but sufficiently close to it that Lemma~\ref{lem:collinear} applies, so
$h^m(G)$ is collinear with $\chain_{q+mR}$ as well as $h^m(V)$. In fact, $G$ is a partition of $K_s$ and we
can use $G$ to code $n$-orbits of $h^m$ unambiguously.

Since by Lemma~\ref{lem:collinear}, $h^m(\chain_q)$ is basically collinear with $\chain_{q+mR}$ (it goes
straight through $\chain_{q+mR}$), and therefore it coils through $\chain_q$ in roughly the same way with the
same number of laps as $T^{mR}_s$ maps $[0, c_1]$ into itself, each $g^i$ intersects $h^m(g^j)$ for at most
$2 \lap(T^{mR}_s) \le 2s^{m(R+\eta)}$ values of $j$. Continuing by induction, we see that each $g^i$
intersects $h^{nm}(g^j)$ for at most $2^n s^{n(mR+\eta)}$ values of $j$, counted `with multiplicity', since
the $h^m$-image of $g^j$ can go through $g^i$ several times. Let us use the partition $(g^i)_{i=0}^{N'}$ to
code the $n$-orbits, \ie for $x \in K_s$, we define an itinerary $e(x) = e_k(x)_{k \ge 0}$ by setting
$e_k(x) = i$ if $h^{mk}(x) \in g^i$. Then there are at most $N' 2^n s^{nm(R+\eta)}$ different itineraries, and
$h^m$-orbits that are $(n, \eps)$-separated must have different length $n$ itineraries. Hence
\begin{eqnarray*}
\htop(h) &=& \frac1m \htop(h^m) \le \frac1m \lim_{n \to \infty} \frac1n
\log N' 2^n s^{nm(R+\eta)} \\
&=& (R+\eta) \log s +  \frac{\log 2}{m} \le R \log s + 2\eta.
\end{eqnarray*}
Here $\eta > 0$ is arbitrary, so the upper bound $\htop(h) \le R \log s$ follows.

Now for the lower bound, let $X_\delta = \{ x  \in [0,s/2] : d(\orb(x), \frac12) > \delta \}$. If
$x \in X_\delta$, then for any $n \in \N$ and maximal interval $J \owns x$ such that $T_s^n|_J$ is monotone,
we have $d(T_s^n(x),\partial T_s^n(J)) > \delta$. It is well-known that $\htop(T_s|_{X_\delta}) \to \log s$
as $\delta \to 0$. Take $\eta > 0$ arbitrary and fix $\delta > 0$ so that
$\htop(T_s|_{X_\delta}) > \log(s-\eta)$. Let $m$ be so large that for any interval $J$ of length
$\ge 2\delta$, $T_s^{mR}(J)$ has at least $(s-2\eta)^{mR}$ branches intersecting $X_\delta$.

Take $\eps  \in (0, 2\delta/s^m)$ and take $q$ so large that $\mesh(\chain_q) < \eps$. For each $q$-point
$x \in C$, let $D_x \subset C$ be the arc of arc-length $2\delta$ centered at $x$, and set
${\mathcal I}_\delta := \bigcup_{x \in E_q} D_x$. Let
$W_1 = h^{-m}(h^m([\bar 0, s_q^1]) \setminus {\mathcal I}_\delta)$; then $W_1$ has at least
$(s-2\eta)^{mR}$ connected components. Continue to define inductively
$W_j = h^{-jm}(h^{jm}(W_{j-1})\setminus {\mathcal I}_\delta)$ for $2 \le j \le n$. The choice of $m$
guarantees that $W_j$ has at least $(s-2\eta)^{jRm}$ connected components. If $x$ and $y$ belong to
different connected components, then there is $0 \le j < n$ and $0 \le k < m$ such that $\pi_0(h^{jm+k}(x))$
and $\pi_0(h^{jm+k}(y))$ are separated by a $\delta$-neighborhood of $\frac12$, so $x$ and $y$ are
$(n, 2\delta/s^m)$-separated for $h^m$. Since $\eps < 2\delta/s^m$, we can select an $(n, \eps)$-separated
set of cardinality $\ge (s-2\eta)^{nRm}$, and therefore
\[
\htop(h) = \frac{1}{m} \htop(h^m) \ge \inf_{\eps > 0} \lim_{n \to \infty}
\frac1{nm} \log
(s-2\eta)^{nRm} = R \log(s-2\eta).
\]
As $\eta$ was arbitrary, the lower bound $\htop(h) \ge R \log s$ follows.
\end{proofof}

We finish with the

\begin{proofof}{Theorem~\ref{thm:entropy_quadratic}}
Let $q_a:x \mapsto 1-ax^2$ be a quadratic map with critical point $0$,
$c_i = q_a^i(0)$ and left fixed point
$u = \frac{1}{2a} (-1-\sqrt{1+4a})$. Assume that $q_a$ is renormalizable,
say it has a period $p$ cycle of
intervals $J_k$, $k = 0, \dots , p-1$ and $u_k \in \partial J_k$ is an
orientation preserving $p$-periodic
point. The structure that emerges from \cite{BaDi3} of the inverse limit
space $\ILq$ of $q_a$ is as follows:
The core of $\ILq$ has $p$ proper subcontinua $G_k$, $k = 0, \dots , p-1$,
which are permuted cyclically by
the shift homeomorphism,  and each of them is homeomorphic to the inverse
limit space of the renormalization $q' := q_a^p|_{J_k}$. One side of the
composant of the point
$(\dots, u_{k-2 \bmod p}, u_{k-1 \bmod p}, u_k)$ compactifies onto the
core of $G_k$, playing the role of the zero composant of 
%$\ILq$ 
$\underleftarrow\lim(J_k,q')$ itself. 
Outside $\cup_k G_k$, the core inverse limit
space has no other subcontinua than points and arcs. Hence, each
homeomorphism $h$ of $\ILq$ can at most
permute the $G_k$ in some way, isotopically to $\sigma^R$, for some $R \in \N_0$. 
However, one can ask whether within $G_k$, $h$ could act
isotopically to $\sigma^{R'}$ where $R' \ne R$ is another(!) power of $\sigma$. 
%For example, if $h$ acts as $\sigma^R$ outside
%$\{ G_k \}$, one might hope that it can still be
%arranged that $h$ maps $G_k$ to $G_{k+R \bmod p}$ isotopically to
In \cite{BS2} we proved that this cannot happen and that every homeomorphism
$h : \ILq \to \ILq$ is isotopic to $\sigma^R$ for some $R \in \Z$.	
\end{proofof}

\medskip
\noindent
Faculty of Mathematics, University of Vienna\\
Oskar Morgensternplatz 1, 1090 Wien, Austria\\
\texttt{henk.bruin@univie.ac.at}\\
\texttt{http://www.mat.univie.ac.at/}$\sim$\texttt{bruin/}\\[3mm]
Department of Mathematics, University of Zagreb\\
Bijeni\v cka 30, 10 000 Zagreb,
Croatia\\
\texttt{sonja@math.hr}\\
\texttt{http://www.math.hr/}$\sim$\texttt{sonja}

\end{document}